\documentclass[reqno,a4wide,11pt]{amsart}
\usepackage[english]{babel}
\usepackage{amsmath}
\usepackage{amssymb}
\usepackage{enumerate}
\usepackage{a4wide}
\usepackage{color}
\usepackage{subfigure}
\usepackage{graphicx}

\newtheorem{thm}{Theorem}[section]

\newtheorem{defin}[thm]{Definition}

\newtheorem{corollary}[thm]{Corollary}

\newtheorem{remark}[thm]{Remark}
\newtheorem{example}[thm]{Example}

\newcommand{\R}{\mathbb{R}}
\newcommand{\Rd}{\mathbb{R}^d}
\newcommand{\dep}{\partial}
\newcommand{\xbar}{\bar{X}}
\newcommand{\dive}{\mathrm{div}}

\newcommand{\sign}{\mathrm{sign}}
\newcommand{\supp}{\mathrm{supp}}
\newcommand{\mespace}{\mathcal{M}_{0}}
\newcommand{\pseudospace}{\mathcal{X}_{0}}
\newcommand{\Xtil}{\widetilde{X}}
\newcommand{\Ftil}{\widetilde{F}}
\newcommand{\rhotil}{\widetilde{\rho}}
\newcommand{\util}{\widetilde{u}}

\newcommand{\be}{\begin{equation}}
\newcommand{\ee}{\end{equation}}
\newcommand{\fer}[1]{\eqref{#1}}
\newcommand{\bee}{\begin{equation*}}
\newcommand{\eee}{\end{equation*}}

\begin{document}

\title{Condensation phenomena in nonlinear drift equations}

\thanks{JAC acknowledges support from the Royal Society by a Wolfson Research
Merit Award and by the Engineering and Physical Sciences Research
Council grant with references EP/K008404/1. JAC was partially
supported by the project MTM2011-27739-C04-02 DGI (Spain) and
2009-SGR-345 from AGAUR-Generalitat de Catalunya. MDF is supported
by the FP7-People Marie Curie CIG (Career Integration Grant)
Diffusive Partial Differential Equations with Nonlocal Interaction
in Biology and Social Sciences (DifNonLoc), by the `Ramon y Cajal'
sub-programme (MICINN-RYC) of the Spanish Ministry of Science and
Innovation, Ref. RYC-2010-06412, and by the by the Ministerio de
Ciencia e Innovaci\'on, grant MTM2011-27739-C04-02.}

\author{Jos\'{e} A. Carrillo}
\address{J. A. Carrillo,
Department of Mathematics,
Imperial College London,
London SW7 2AZ}
\email{carrillo@imperial.ac.uk}

\author{Marco Di Francesco}
\address{Marco Di Francesco, University of Bath. Mathematical Sciences, University of Bath, Claverton Down, Bath (UK), BA2 7AY.}
\email{m.difrancesco@bath.ac.uk}

\author{Giuseppe Toscani}
\address{Giuseppe Toscani, Dipartimento di Matematica, Universit\`{a} di Pavia, Via Ferrata 1, 27100 Pavia, Italy}
\email{giuseppe.toscani@unipv.it}

\smallskip

\begin{abstract}
We study nonnegative, measure-valued solutions to nonlinear drift
type equations modelling concentration phenomena related to
Bose-Einstein particles. In one spatial dimension, we prove
existence and uniqueness for measure solutions. Moreover, we prove
that all solutions blow up in finite time leading to a
concentration of mass only at the origin, and the concentrated
mass absorbs increasingly the mass converging to the total mass as
$t\to\infty$. Our analysis makes a substantial use of independent
variable scalings and pseudo-inverse functions techniques.
\end{abstract}



\maketitle

\section{Introduction}

In this paper we perform a rigorous study of nonnegative measure solutions
of the nonlinear drift-type equation
\begin{equation}\label{eq:main}
    \dep_\tau f = \nabla_v\cdot( v f(1+f^\gamma)) ,\qquad \gamma > 0,
\end{equation}
posed on $v\in \R^d$, $d\geq 1$, and $\tau\geq 0$, with initial condition $f_I\in
L^1_+(\R^d)$.

In the case $\gamma =1$, the nonlinear drift on the right-hand side of
equation \fer{eq:main} appears as the confinement operator in the
Fokker-Planck equation introduced by Kaniadakis and Quarati
\cite{KQ93,KQ94}, who proposed a correction to the linear drift term of
the Fokker-Planck equation to account for the presence of quantum
indistinguishable particles, bosons or fermions. There, the evolution
equation for bosons density has been postulated in the form
 \begin{equation}
\label{FPB}
   \dep_\tau f=  \nabla\cdot \left[ \nabla f + vf(1+\epsilon
f)\right],\qquad \epsilon>0.
 \end{equation}
The fundamental assumption leading to the correction in the
drift term of the linear Fokker-Planck equation is established in
\cite{CC70} by imposing that a state of congestion emerges when the mean distance between
neighboring particles is comparable to the size of the quantum wave
fields in which the particles are embedded. In the case of a gas composed by Bose-Einstein identical particles,
(according to quantum theory) the probability that a new particle will enter the velocity range $dv$ increases when the neighboring particles are of the same type. More precisely, the presence of $f(v) dv$ particles per unit
volume increases this probability in the ratio $1+\epsilon f(v)$, for some $\epsilon>0$.

By a direct calculuation, one can easily check that the Bose-Einstein distribution \cite{CC70}
\be\label{BEsteady}
 f_{\lambda} (v) = \frac 1\epsilon \left[ e^{|v|^2/2 +\lambda} - 1 \right]^{-1},
\ee
satisfies the equation
 \[
\nabla f_\lambda(v) + vf_\lambda(v)(1+\epsilon f_\lambda(v)) = 0
 \]
for any fixed positive constant $\lambda$, and it is therefore a stationary state for the equation
\fer{FPB}. The constant $\lambda$ is
related to the mass of Bose-Einstein distribution
 \[
m_\lambda = \int_{\R^d}\frac 1\epsilon \left[ e^{|v|^2/2 +\lambda} - 1
\right]^{-1} \, dv,
 \]
and, since $m_\lambda$ is decreasing as $\lambda$ increases, the
maximum value of $m_\lambda$ is attained at $\lambda =0$. If $d
=3$, the value
$$ 
m_c = m_0 = \int_{\R^3}\frac
1\epsilon \left[ e^{|v|^2/2} - 1 \right]^{-1} \, dv < +\infty
$$
defines a \emph{critical mass}.

One of the main problems in kinetic equations
relaxing towards a stationary state characterized by the existence of a
critical mass is to detect a singular
part (the condensate) for the solution when the initial distribution
has a supercritical mass $m> m_c$,. In general this phenomenon is heavily
dependent of the dimension of the physical space. For instance, in dimension $d
\le 2$ the maximal mass $m_0$ of the Bose-Einstein
distribution \fer{BEsteady} is $+\infty$, and the eventual formation
of a condensate is lost.

The kinetics of Bose--Einstein condensation - in particular, the way in which the
Bose fluid undergoes a transition from a normal fluid to one with a
condensate component - has been the object of various investigations
\cite{EM99,EM01,EMV03,JPR06,LLPR,ST95,ST97,Spo08}. The results in the aforementioned literature are
mainly based on study of kinetic equations, like the quantum Boltzmann
equation and the Boltzmann--Nordheim equation, which describes the
dynamics of weakly interacting quantum fluids, see also the recent \cite{EscVel}.

Several kinetic models for Bose-Einstein particles have been proposed and studied
in the literature. In particular, a related model described by means of
Fokker-Planck type non linear operators has been proposed by
Kompaneets \cite{Kom57} to describe the evolution of the radiation
distribution in a homogeneous plasma when radiation interacts with
matter via Compton scattering. This equation has been exhaustively
studied in \cite{EHV98}. It should be noticed that the Kompaneets
equation also features a nonlinear drift of similar type as in \fer{eq:main}.

The study of the mathematical theory for the Fokker-Planck equation \fer{FPB}, both in the
bosons and fermions cases, is quite recent, and still partly open. Existence of solutions for
the fermions case has been investigated in \cite{CLR08}. For bosons,
equation \fer{FPB} has been studied in dimension one of the velocity
variable  in \cite{CRS08}. In this case, however, the equilibrium
Bose-Einstein density is a smooth function, which makes it possible to
prove exponential convergence to equilibrium by relying on standard entropy
methods. The three-dimensional problem, and the eventual formation of a
condensed part in the solution of the equation \fer{FPB}, has been recently
investigated in \cite{To12}, where it is shown that, for a given initial
supercritical mass, there is blow-up in finite time of the solution.
However, nothing has been proven on the behavior of the solution after the
blow-up time. It is clear that both the finite-time blow up proven in
\cite{To12}, and the possible subsequent formation of a condensed part,
are due to the presence of the nonlinearity in the drift equation. Indeed,
the dissipation of energy in the linear drift equation is not enough to
produce blow up in finite time, and concentration of mass is obtained only
as time goes to infinity.

While the case $\gamma = 1$ is the most relevant one from the
physical point of view (due to its relationship with the
Bose-Einstein distribution), Fokker-Planck equations of type \fer{eq:main},
with $\gamma >1$ and with linear diffusion have
been considered in \cite{BGT11}, with the aim of finding minimizers to
sub-linear entropies. There, it has been shown that a suitable
coupling of the degree of nonlinearity of the drift and of the space
dimension ($\gamma >2$ if $d=1$) allows the equilibrium profile to
have a singular part, which is localized in space in a precise way.

The analysis of \cite{BGT11,To12} motivates the study of equation
\fer{eq:main}, both in dimension one, and in higher dimensions of
the velocity variables. By enlightening its main properties - in particular towards the direction of the formation of condensation in finite time, and the continuation of measure solutions for all times - we hope to shed a light on
the main problem of the formation of a condensed part in the case
of the complete Fokker-Planck equation with linear diffusion \fer{FPB}.

The first step in our analysis is to reduce \eqref{eq:main} to the new equation
\begin{equation}\label{eq:new_intro}
    \dep_t \rho - \dive (x\rho^{1+\gamma}) = 0,
\end{equation}
via a time dependent scaling which removes the linear drift part,
see section \ref{sec:multi_d} below. The main advantage of
\eqref{eq:new_intro} is that the nonlinear drift term becomes
\emph{homogeneous}.

In one space dimension, this very simple key observation paves the
way to our next key remark, namely that the Cauchy problem for the equation
\eqref{eq:new_intro} can be related to the Cauchy-Dirichlet
problems for a classical scalar conservation law
\begin{equation}\label{eq:burgers_intro}
\begin{cases}
    u_t + \left(\frac{1}{1+\gamma} u^{1+\gamma}\right)_\xi = 0& \hbox{if}\quad \xi<0\\
    u_t - \left(\frac{1}{1+\gamma} u^{1+\gamma}\right)_\xi = 0& \hbox{if}\quad \xi>0,
\end{cases}
\end{equation}
with boundary datum $u(0,t)=0$, via the scaling
\begin{align}
    & u(\xi,t)=x'(\xi) \rho\left( x(\xi),t\right),\quad \xi\in \R,\quad t\geq 0\nonumber\\
  & x(\xi)= \sign(\xi)\frac{1}{1+\gamma}(\gamma
    |\xi|)^{(1+\gamma)/\gamma},\qquad x'(\xi)=(\gamma |\xi|)^{1/\gamma}.\label{eq:scaling_intro}
\end{align}
see section \ref{sec:oned} below.

The scaling \eqref{eq:scaling_intro} can be easily justified in
case of a bounded solution $\rho\in L^\infty$ of
\eqref{eq:new_intro}, which corresponds to $u(0,t)=0$. Then, when the
left and right traces of $u$ at $\xi=0$ become positive, $u$
and $\rho$ start losing mass. Based on the previous results in the
literature, and having not included diffusion in the model, we
clearly expect that the lost mass gets concentrated. Hence, the
scaling \eqref{eq:scaling_intro} becomes the main tool to explore
condensation phenomena for \eqref{eq:main} via the auxiliary
equation \eqref{eq:new_intro}.

On the other hand, this requires a global-in-time well-posedness
theory for \emph{measure solutions} for \eqref{eq:main}, in a way
to embrace also entropy solutions a la Oleinik--Kru{\v{z}}kov
\cite{oleinik63,kruzkov}, at least for short times. The
achievement of such a goal constitutes the main result of this
paper, and it is rigorously stated in Definition
\ref{def:entropy_measure} and in Theorem \ref{thm:main} in Section
\ref{sec:existence}, the main difficulty being represented by the
nonlinearity $\rho^{1+\gamma}$ in the drift term, which renders
the definition of a Dirac delta type solution a non trivial task.
The main tool used to overcome this difficulty is to consider the
\emph{pseudo--inverse} formulation of the equation
\eqref{eq:new_intro}, see section \ref{sec:measure}, which allows
to consider constant solutions in the pseudo--inverse variable
corresponding to concentrated parts in the solution to
\eqref{eq:new_intro}.

Based on our existence and uniqueness result in Theorem
\ref{thm:main}, we then analyse the large time behaviour or the
measure solutions. More precisely, under suitable restriction on
the initial conditions (still in a context of \emph{large} data)
we prove that all measure solutions to \eqref{eq:new_intro}
concentrate in a finite time, and the concentrated mass is
strictly increasing after the blow-up time. Finally, we prove that
\emph{all} the total mass concentrates as $t\rightarrow +\infty$.
This is proven in Theorem \ref{thm:blow_up}. A quite enlightening
special solution is described in the example \ref{exe:explicit},
in which we see a bounded initial data evolving towards a
condensate state as follows: at some time $t^*$ the profile of the
solution \emph{becomes unbounded at zero}, but still integrable.
Immediately after $t^*$, the solution concentrates gradually, and
all the mass is concentrated in infinite time.

It is interesting to notice that (partially) explicit solutions of
the pseudo-inverse equation \eqref{eq:pseudo_scaled_measures} are
provided in section \ref{sec:measure} by the method of
characteristics for fully nonlinear PDEs, in the spirit of
Hamilton-Jacobi type equations, see \cite{lions}. A possible way
to achieve existence and uniqueness for that problem was then the
use of the theory of viscosity solutions, see e. g.
\cite{crandall}. On the other hand, the special structure of the
equation and its relation to the scalar conservation law
\eqref{eq:burgers_intro} allow to prove existence and uniqueness
in a more direct way, and this gives as a by product the
equivalence with Oleinik's entropy solutions of
\eqref{eq:new_intro} for short times.

The paper is organized as follows. In section \ref{sec:multi_d} we
present the scaling needed to write \eqref{eq:new_intro}, and
provide a short time existence result of smooth solutions via
characteristics, which holds in arbitrary space dimension. In section \ref{sec:oned} we analyse the
one-dimensional case, and explain in detail the relation with the
classical conservation law model \eqref{eq:burgers_intro}. Here we
provide also the special solution mentioned above. In section
\ref{sec:measure} we introduce the use of the pseudo-inverse
variable, and provide explicit solutions via the method of
characteristics for fully nonlinear PDEs. In section
\ref{sec:existence} we provide out main results in Theorem
\ref{thm:main}. In section \ref{sec:qualitative} we prove the results
about the qualitative and asymptotic behaviour mentioned before.

Let us mention that all the main results are proven for the
one-dimensional scaled equation \eqref{eq:new_intro}, but they can
be easily re-stated for the solution $f$ to \eqref{eq:main}, see
Remarks \ref{rem:original1} and \ref{rem:original2}.

\section{The multi dimensional case}\label{sec:multi_d}

Throughout the whole paper, $\mathcal{P}(\R^d)$ will denote the
space of probability measures on $\R^d$. $\mathcal{L}^d$ denotes
the Lebesgue measure on $\R^d$. $\delta_{x_0}$ is the usual Dirac
delta measure centered at $x_0$. For a given function $\rho \in
L^1(\Rd)$ we use the notation
$$
\supp(\rho):=\left\{x\in \R^d:\
\int_{B(x,\epsilon)}\rho(x)\,dx\neq 0 \mbox{ for all }
\epsilon>0\right\}\,.
$$

We start with the observation that there exists a time-dependent
mass--preserving scaling which removes the linear drift term
$\partial_v(vf)$ from \eqref{eq:main}. More precisely, let us set
\begin{equation}\label{eq:scaling}
\begin{cases}
    f(v,\tau)=e^{d\tau}\rho( x,t) & \\
    x =e^\tau v & \\
    t = \frac{e^{d\gamma\tau} - 1}{d\gamma}.
\end{cases}
\end{equation}
We obtain the following scaled equation for $\rho(x,t)$
\begin{equation}\label{eq:main_scaled}
    \dep_t \rho - \dive (x\rho^{1+\gamma}) = 0
\end{equation}
with $\rho(x,0)=f_I(x)$. The new equation \eqref{eq:main_scaled}
has the great advantage of featuring a \emph{homogeneous
nonlinearity in the drift term}. We next write
\eqref{eq:main_scaled} in non conservative form
\begin{equation}\label{eq:main_scaled_noncon}
    \rho_t -(1+\gamma)\rho^{\gamma} x \cdot \nabla \rho = d \rho^{1+\gamma},
\end{equation}
which can only be done in case of smooth solutions. For a fixed
$x_0\in \R^d$, the characteristics curves for the nonlinear
first-order equation \eqref{eq:main_scaled_noncon} starting at
$x_0$ are denoted by $X_{x_0}(t)$, with
$U_{x_0}(t):=\rho(X_{x_0}(t),t)$. They solve
\begin{equation}\label{eq:charact1}
    \begin{cases}
    \dot{X}_{x_0}(t) = -(1+\gamma) X_{x_0}(t) U_{x_0}(t)^{\gamma}, & X_{x_0}(0)= x_0
    \\[2mm]
    \dot{U}_{x_0}(t) = d U_{x_0}^{1+\gamma}(t), & U_{x_0}(0)=u_0:=f_I(x_0).
    \end{cases}
\end{equation}
The solution to the second equation in \eqref{eq:charact1} is
\begin{equation*}
    U_{x_0}(t)= \frac{f_I(x_0)}{(1- \gamma d f_I^\gamma(x_0) t)^{1/\gamma}},
\end{equation*}
which implies
\begin{equation}\label{eq:char_multid}
\begin{cases}
     X_{x_0}(t) = x_0 (1-\gamma d f_I^\gamma(x_0) t)^{(1+\gamma)/\gamma d}&
     \\[3mm]
     \rho(X_{x_0}(t),t)= \frac{f_I(x_0)}{(1- \gamma d f_I^\gamma(x_0) t)^{1/\gamma}}.&
\end{cases}
\end{equation}
The above formulae \eqref{eq:char_multid} suggest two important
properties to hold for \eqref{eq:main_scaled} in any dimension:
\begin{enumerate}
  \item[i)] all solutions become unbounded in a finite time $t^*$ depending on
the initial condition,
  \item[ii)] the only point $x$ at which a solution may
become unbounded is $x=0$.
\end{enumerate}
Another property which is suggested by \eqref{eq:char_multid} is the \emph{confinement} property
\begin{equation*}
  \supp[\rho(\cdot,t)]\subset \mathrm{Conv}\left(\supp(f_I)\cup \{0\}\right) \qquad \hbox{for all}\ t>0,
\end{equation*}
which formally holds as long as solutions can be obtained via
characteristics, since all the curves $X_{x_0}(\cdot)$ \emph{point
towards} the axis $x=0$. Here, $\mathrm{Conv}(A)$ denotes the
convex hull of $A$. Finally, \eqref{eq:char_multid} suggests that
$X_0(t)\equiv 0$ as long as the solution $\rho$ stays smooth.

As usual in the context of scalar conservation laws
\cite{dafermos}, we can use of the characteristic curves recovered
in \eqref{eq:char_multid} to produce a local existence theorem for
smooth solutions to \eqref{eq:main_scaled} with smooth initial
data. To perform this task, one has to invert the formula
\eqref{eq:char_multid} for characteristics by recovering the
initial point $x_0$ as a function of $X_{x_0}(t)$ and $t$. By the
implicit function theorem, this can be done as long as the
function
\begin{equation*}
    \mathcal{F}(x,x_0,t):=x_0 (1-\gamma df_I^\gamma(x_0)t)^{(1+\gamma)/\gamma d} - x
\end{equation*}
satisfies
\begin{equation*}
    \frac{\mathrm{D} \mathcal{F}}{\mathrm{D} x_0}\neq 0.
\end{equation*}
A straightforward computation yields
\begin{equation}\label{eq:dini}
   \frac{\mathrm{D} \mathcal{F}}{\mathrm{D} x_0} = (1-\gamma d f_I^\gamma(x_0) t)^{(1+\gamma)/\gamma d}\mathbb{I}_d
    - (1+\gamma)x_0 \otimes \nabla f_I^\gamma (x_0) (1-\gamma d f_I^\gamma (x_0) t)^{(1+\gamma-\gamma d)/\gamma d} t.
\end{equation}
Since $\frac{\dep \mathcal{F}}{\dep x_0}|_{t=0}= \mathbb{I}_d$, it
is clear that characteristics do not cross at least on a small
time--strip provided $\nabla \rho^\gamma_0$ is finite almost
everywhere. In one space dimension, the above condition
\eqref{eq:dini} reads
\begin{equation*}
    (1-\gamma f_I^\gamma(x_0)t)^{(1+\gamma)/\gamma} -(1+\gamma)x_0(1-\gamma f_I^\gamma(x_0)t)^{1/\gamma}(f_I^\gamma)'(x_0) t \neq 0,
\end{equation*}
which is satisfied for all $0\leq t < (\gamma
f_I^\gamma(x_0))^{-1}$ provided $x_0f_I'(x_0) \leq 0$. Therefore,
the latter condition ensures that the solution $\rho$ is classical
on $(x,t)\in \R \times [0,t^*)$ with $t^*=(\gamma \max
f_I^\gamma)^{-1}$. At $t=t^*$, $\rho$ tends to $+\infty$ at the
point $x=0$.

In more than one space dimension, and for a radially symmetric
initial datum $f_I(x_0)=\widetilde{f}_I(|x_0|)$, \eqref{eq:dini}
reads
\begin{equation*}
   \frac{\mathrm{D} \mathcal{F}}{\mathrm{D} x_0} = (1-\gamma d \widetilde{f}_I^\gamma(|x_0|) t)^{(1+\gamma)/\gamma d}\mathbb{I}_d
   - (1+\gamma)x_0 \otimes x_0 \frac{\left(\widetilde{f}_I^\gamma\right)_r (|x_0|)}{|x_0|}
   (1-\gamma d \widetilde{f}_I^\gamma (|x_0|) t)^{(1+\gamma-\gamma d)/\gamma d} t.
\end{equation*}
Therefore, in case $\rho$ is initially radially non-increasing,
the solution is classical until $t=(\gamma d \max
f_I^\gamma)^{-1}$. This is due to the fact that the matrix
$x_0\otimes x_0$ is non-negative definite. We have therefore
proved the following

\begin{thm}[Local existence of smooth solutions]\label{thm:local}
Let $f_I \in C^1\cap L^\infty(\Rd)$, $f_I \geq 0$. Then, there
exist a time $t^*>0$ and $\rho(\cdot,\cdot)\in C^1([0,t^*)\times
\Rd)$ solution to \eqref{eq:main_scaled_noncon}. Moreover,
\begin{itemize}
  \item if $f_I$ is radially non increasing, then the maximal time $t^*$
is given by
\begin{equation}\label{eq:maximal_time}
    t^* =(\gamma d \max f_I^\gamma)^{-1},
\end{equation}
    \item if the maximal time is given by \eqref{eq:maximal_time}, then the solution $\rho$ blows up at $x=0$, namely
\begin{equation*}
    \lim_{(x,t)\rightarrow (0,t^*)}\rho(x,t)=+\infty,
\end{equation*}
  \item if the maximal time satisfies $t^*<(\gamma d \max
  f_I^\gamma)^{-1}$, then the solution $\rho$ develops a
  discontinuity before blowing up.
\end{itemize}
\end{thm}

The result in Theorem \ref{thm:local} shows that discontinuities
may arise in a finite time, as usual in this context, since
characteristics may cross in finite time. The result in Theorem
\ref{thm:local} can be easily re-formulated in terms of the
original equation \eqref{eq:main}. We omit the details.

\section{The one dimensional case}\label{sec:oned}

Let us then recover a reasonable notion of weak solution at least
in the one dimensional case. Being the only possible concentration
point at zero, the characteristics all pointing towards the
origin, and dealing with a conservation law, we cannot expect
anything else than a concentration of mass at zero. Then, we need
a notion of measure solution that allows for this possibility. We
are going to start by introducing a notion of weak solution that
avoids the concentration at the origin issue by allowing a loss of
mass through the origin. This will be just an intermediate step in
our construction.

\begin{defin}[Weak solutions]\label{def:weak_sol}
Let $T>0$. Then, $\rho(x,t)$ is a weak solution to
\eqref{eq:main_scaled} on $[0,T]$ with initial datum $f_I \in
L^\infty(\R)$ if and only if
\begin{itemize}
  \item $x\rho^{1+\gamma}\in L^\infty(\R \times [0,T))$
  \item for all $\varphi\in C_c^\infty(\R \times [0,T))$ such that $\varphi(0)=0$, one has
  \begin{align}
    & \int_{\R^d}\int_0^{T} \varphi_t(x,t) \rho(x,t) dx dt = \int_{\R^d}\varphi(x,0)f_I(x) dx + \int_{\R^d}\int_0^{T} \rho^{1+\gamma}
    (x,t) x\cdot \nabla \varphi (x,t) dx dt .\label{eq:def_weak}
  \end{align}
\end{itemize}
\end{defin}

We will clarify in Section \ref{sec:oned} that examples of non
uniqueness of weak solutions can be easily found. Therefore a
definition of entropy solution in the spirit of \cite{kruzkov} is
needed. This notion can be written for the evolution equation
\eqref{eq:main_scaled} but we prefer to postpone it since it will
be much clearer through a new change of variables.

In order to deal with entropy solutions, we shall restrict to the
case of initial condition $f_I\in BV(\R)$. Note that in one
dimension, this implies that $f_I\in L^1\cap L^\infty(\R)$. We
will assume that we work with probability measures denoted by
$\mathcal{P}(\R)$, and thus with unit mass non-negative densities.
Further restrictions will be required later on.

Let $f_I\in \mathcal{P}(\R)\cap BV(\R)$ with the additional
condition $f_I\in C^1(\R)$, and let $\rho$ be the local-in-time
$C^1$ bounded solution to
\begin{equation}\label{eq:main_oned}
    \rho_t - ( x \rho^{1+\gamma})_x =0
\end{equation}
with initial datum $f_I$ provided by Theorem \ref{thm:local}. As
we already pointed out in section \ref{sec:multi_d}, the
characteristic curve generated at $x_0=0$ is a vertical line, and
it therefore separates the mass on $x<0$ from the mass on $x>0$.
More precisely, let $t^*$ be the maximal time of existence of $\rho$
provided in Theorem \ref{thm:local}, for all $t\in [0,t^*[$ we
have
\begin{equation*}
  \int_{-\infty}^0 \rho(x,t) dx = m_L:= \int_{-\infty}^0 f_I(x) dx,\qquad \int_0^{+\infty} \rho(x,t) dx = m_R:= \int_0^{+\infty} f_I(x) dx
\end{equation*}
This property is due to the expression for the flux $x\rho^{1+\gamma}$, which clearly vanishes at $x=0$ as long as $\rho$ is bounded.
We introduce the scaled solution $u$ as follows:
\begin{align}
    & u(\xi,t)=x'(\xi) \rho\left( x(\xi),t\right),\quad \xi\in \R,\quad t\geq 0\nonumber\\
  & x(\xi)= \sign(\xi)\frac{1}{1+\gamma}(\gamma
    |\xi|)^{(1+\gamma)/\gamma},\qquad x'(\xi)=(\gamma |\xi|)^{1/\gamma}.\label{eq:change_inverse}
\end{align}
It is immediately checked that
\begin{align*}
    & \int_0^{+\infty}\rho(x,t) dx = \int_0^{+\infty}
    u(\xi,t)d\xi,\\
    & \int_{-\infty}^{0}\rho(x,t) dx =
    \int_{-\infty}^{0}  u(\xi,t)d \xi,
\end{align*}
and that $u$ satisfies
\begin{equation}\label{eq:burgers0}
    u_t - \sign(\xi)\left(\frac{1}{1+\gamma} u^{1+\gamma}\right)_\xi = 0.
\end{equation}
Since $\rho$ is bounded, $u$ has the trace property
\begin{equation}\label{eq:u_zero_trace}
  u(0,t)=0\quad \hbox{for all}\ t\in [0,t^*[.
\end{equation}
The change of variable $x=x(\xi)$ introduced in \eqref{eq:change_inverse} admits the inverse transformation
\begin{equation}\label{eq:change_x_to_xi}
  \xi(x)=\sign(x) \frac{1}{\gamma}[(1+\gamma)|x|]^{\frac{\gamma}{1+\gamma}},
\end{equation}
which is not differentiable at $x=0$. However, in view of \eqref{eq:u_zero_trace}, the inverse change of variable
\begin{equation}\label{eq:change_to_burgers}
    \rho(x,t)= \xi'(x) u(\xi(x),t)
\end{equation}
is well defined for all $x\neq 0$.

Although the change of variable \eqref{eq:change_inverse} has the
advantage of being defined on $\xi\in \R$, the equation
\eqref{eq:burgers0} may be uncomfortable to deal with because of
the discontinuous coefficient at $\xi=0$. We shall therefore
consider two separate Cauchy problems:
\begin{align}
    &
    \begin{cases}
    u_t - \left(\frac{1}{1+\gamma} u^{1+\gamma}\right)_\xi = 0 &
    \xi>0,\ \ t\geq 0 \\
    u(\xi,0)=u_{I,R}(\xi):=(\gamma\xi)^{1/\gamma}f_I\left(\frac{(\gamma
    \xi)^{(1+\gamma)/\gamma}}{1+\gamma}\right) & \xi>0,
    \end{cases}\label{eq:IBV1}\\
&
    \begin{cases}
    u_t + \left(\frac{1}{1+\gamma} u^{1+\gamma}\right)_\xi = 0 &
    \xi<0,\ \ t\geq 0 \\
    u(\xi,0)=u_{I,L}(\xi):=(-\gamma\xi)^{1/\gamma}f_I\left(-\frac{(-\gamma
    \xi)^{(1+\gamma)/\gamma}}{1+\gamma}\right) & \xi<0,
    \end{cases}\label{eq:IBV2}
\end{align}

We shall denote solutions to \eqref{eq:IBV1} (to \eqref{eq:IBV2} resp.) by $u_R$ ($u_L$ resp.), and we shall often refer to $u$ as
\begin{equation}\label{eq:u_def_split}
  u(\xi,t)=
\begin{cases}
u_L(\xi,t) & \xi<0\\
u_R(\xi,t) & \xi>0
\end{cases}
\end{equation}

As the initial condition in both problems \eqref{eq:IBV1} and
\eqref{eq:IBV2} is nonnegative, the characteristic curves
generated near the boundary $\xi=0$ are \emph{outgoing}.
Therefore, no boundary condition needs to be prescribed in order
to achieve a unique entropy solution for the initial value
problems \eqref{eq:IBV1}-\eqref{eq:IBV2}, see
\cite{bardos_et_al,dubois}. Moreover, it is well known in the
context of scalar conservation laws with uniformly convex or
concave flux that the classical notion of entropy solution in
Kru{\v{z}}kov sense \cite{kruzkov} is equivalent to the notion of
weak solutions together with the \emph{Oleinik's condition}
\cite{oleinik63}. This condition roughly speaking requires a
one-sided bound for the space derivative in a distributional
sense. In the present context, since the convexity of the flux
depends on the sign of $\xi$, and in view of well known results
cf. \cite{hoff83}, the correct Oleinik type condition is given by
the \emph{distributional} inequality
\begin{equation*}
    -\sign(\xi)(u^{\gamma})_{\xi}\leq \frac{C}{t},\quad \xi\neq 0,\ \
    t>0.
\end{equation*}
for some positive constant $C$. However, our assumption of $f_I\in BV(\R)$ allows to formulate the entropy condition in the even simpler way
\begin{align}
    & \lim_{\xi\nearrow \xi_0^-}u(\xi,t) \geq \lim_{\xi\searrow \xi_0^+}u(\xi,t)\qquad \hbox{if}\: \xi_0<0,\nonumber\\
     & \lim_{\xi\nearrow \xi_0^-}u(\xi,t) \leq \lim_{\xi\searrow \xi_0^+}u(\xi,t)\qquad \hbox{if}\: \xi_0>0.
    \label{eq:oleinik_simple}
\end{align}

\begin{defin}\label{def:entropy_burgers}
Let $f_I\in \mathcal{P}(\R)\cap BV(\R)$. A function $u\in
L^\infty([0,+\infty);L^1(\R)\cap BV(\R))$ is an
\emph{entropy solution} to the Cauchy problem
\eqref{eq:IBV1}--\eqref{eq:IBV2} if $u$ can be written as
\eqref{eq:u_def_split} with $u_R$ solving \eqref{eq:IBV1} and
$u_L$ solving \eqref{eq:IBV2} in a weak sense, and $u$ satisfies
the jump conditions \eqref{eq:oleinik_simple}.
\end{defin}

 We can therefore collect the results
in \cite{oleinik63,kruzkov,bardos_et_al,hoff83} adapted to problem
\eqref{eq:IBV1}-\eqref{eq:IBV2} in the following theorem.

\begin{thm}\label{thm:burgers}
Let $f_I\in \mathcal{P}(\R)\cap BV(\R)$, then there exists a
unique $u\in L^\infty([0,+\infty);L^1(\R)\cap BV(\R))$
entropy solution to the problem \eqref{eq:IBV1}--\eqref{eq:IBV2}
in the sense of Definition \ref{def:entropy_burgers}.
\end{thm}

\begin{remark}[Comparison principle]\label{rem:comparison}
\emph{As a consequence of the classical theory of scalar
conservation laws (see \cite[Chapter 6]{dafermos}), for two given
entropy solutions $u_L^1,u_L^2$ to \eqref{eq:IBV1} we have
$u_L^1(x,t)\leq u_L^2(x,t)$ for almost every $x<0$ and for all
$t>0$ provided $u_L^1(x,0)\leq u_L^2(x,0)$.}
\end{remark}

\begin{remark}[Mass conservation for solutions to \eqref{eq:main_scaled}]\label{rem:zerotrace}
\emph{It can be easily checked (see the proof of Theorem
\ref{thm:main} below) that if $u$ is a entropy solution to
\eqref{eq:IBV1}--\eqref{eq:IBV2} then $\rho$ defined from $u$ via
the scaling \eqref{eq:change_to_burgers} is a weak solution to
\eqref{eq:main_scaled} on $[0,+\infty]$ according to Definition
\ref{def:weak_sol}.}

\emph{If the two traces $u(0^-,t)$ and $u(0^+,t)$ of the entropy
solutions to \eqref{eq:IBV1}--\eqref{eq:IBV2} are zero, the total
mass of the system is conserved. However, this property is lost as
soon as one of the two traces becomes positive. In this case,
$\rho$ solution to \eqref{eq:main_oned} becomes unbounded at $x=0$
due to the change of variables \eqref{eq:change_to_burgers} and
start losing mass at the origin. Our conjecture is that, speaking
in terms of $\rho$, such lost mass gets concentrated to a delta
measure at the origin, as we shall see in the example below.
Therefore, the notion of solution $u$ provided in definition
\ref{def:entropy_burgers} will be used to extend our solutions
$\rho$ in a measure sense after their blow up in a finite time.}
\end{remark}

\begin{example}[An explicit solution]\label{exe:explicit}
\emph{We shall provide an
explicit solution to \eqref{eq:main_oned} by setting as initial datum
\begin{equation}\label{eq:initial_datum_explicit}
    \rho(x,0)=f_I(x)=
    \begin{cases}
    1 & \hbox{if}\ \ 0\leq x \leq \frac{1}{1+\gamma} \\
    0 & \hbox{otherwise}.
    \end{cases}
\end{equation}
We use the scaling \eqref{eq:change_inverse}, which turns the initial condition $f_I$ into
\begin{equation*}
  u_0(\xi):=(\gamma \xi)^{1/\gamma}f_I\left(\frac{1}{1+\gamma}(\gamma
    |\xi|)^{(1+\gamma)/\gamma}\right) = (\gamma \xi)^{1/\gamma}\chi_{[0,1/\gamma]}
\end{equation*}
We recall the characteristic equations for the scaled equation \eqref{eq:IBV1}:
\begin{equation*}
  \left\{
    \begin{array}{ll}
      \xi_{\xi_0}(t)= \xi_0 (1- \gamma t), &  \\
      u(\xi_{\xi_0}(t),t)=u_0(\xi_0), &
    \end{array}
  \right.
\end{equation*}
for $0\leq \xi_0\leq 1/\gamma$. If $\xi_0 >1/\gamma$ we have
$\xi_{\xi_0}(t)\equiv \xi_0$. There is a rarefaction fan $1/\gamma
-t <\xi<1/\gamma$ which is filled up with the characteristic lines
$\xi_\lambda(t)=\frac{1}{\gamma}-\lambda t$, with $\lambda\in
(0,1)$. The profile in the rarefaction fan is
$$
u(\xi,t)=\left(\frac{1-\gamma \xi}{\gamma t}\right)^{1/\gamma}.
$$
Therefore, the unique entropy solution to \eqref{eq:IBV1} with
initial datum $u_0$ is
\begin{equation*}
u(\xi,t)=
  \left\{
     \begin{array}{ll}
       \left(\frac{\gamma \xi}{1-\gamma t}\right)^{1/\gamma}, & \hbox{if}\ 0\leq \xi \leq \frac{1}{\gamma}-t,\ \ \hbox{and}\ 0\leq t <\frac{1}{\gamma} \\
       \left(\frac{1-\gamma \xi}{\gamma t}\right)^{1/\gamma}, & \hbox{if}\ \max\{0,\frac{1}{\gamma}-t\} \leq \xi \leq \frac{1}{\gamma} \\
       0, & \hbox{otherwise}
     \end{array}
   \right.
\end{equation*}
By applying the scaling \eqref{eq:change_to_burgers}, we obtain
the solution $\rho$ to \eqref{eq:main_oned} with initial datum
$f_I$ given in \eqref{eq:initial_datum_explicit}:
\begin{equation}\label{eq:solution_explicit}
  \rho(x,t)=
\left\{
  \begin{array}{ll}
    \left(\frac{1}{1-\gamma t}\right)^{\frac{1}{\gamma}}, & \hbox{if}\ \ 0\leq x \leq \frac{(1-\gamma t)^{\frac{\gamma}{1+\gamma}}}{1+\gamma}\ \ \hbox{and}\ \ 0\leq t \leq 1/\gamma \\
    \left[\frac{((1+\gamma)x)^{-\frac{\gamma}{1+\gamma}}-1}{\gamma t}\right]^{\frac{1}{\gamma}}, & \hbox{if}\ \ \max\{0,\frac{(1-\gamma t)^{\frac{\gamma}{1+\gamma}}}{1+\gamma}\}\leq x \leq \frac{1}{1+\gamma} \\
    0, & \hbox{otherwise}
  \end{array}
\right.
\end{equation}
It is easily seen that $\rho$ in \eqref{eq:solution_explicit}
preserves the total mass $\int_0^{+\infty}\rho(x,t)dx =
\frac{1}{1+\gamma}$ as long as $t\leq \frac{1}{\gamma}$. At
$t=t^*:=\frac{1}{\gamma}$, the total mass starts decreasing. More
precisely, for $t\geq t^*$,
\begin{equation*}
  m(t):=\int_0^{+\infty}\rho(x,t) dx = \left(\frac{1}{\gamma t}\right)^{1/\gamma}\frac{1}{1+\gamma}.
\end{equation*}
Such a behaviour suggests that the mass of $\rho(\cdot,t)$ starts
concentrating at $t>t^*$, in particular that a Dirac's delta part
centered at zero with mass
$$
\frac{1}{1+\gamma}\left(1-\left(\frac{1}{\gamma
t}\right)^{1/\gamma}\right)
$$
is generated for $t>\frac{1}{\gamma}$. Please notice that the mass
concentrated at zero is continuous in time, and it converges to
the total mass of $\rho$ as $t\rightarrow +\infty$. The solution
found here can be easily scaled back to a solution $f$ to
\eqref{eq:main} via the scaling \eqref{eq:scaling}. The
qualitative behaviour is substantially the same, with the only
difference that the $f$ is subject to a confinement, which makes
the support of $f$ shrink to $v=0$ for large times with algebraic
rate.}
\end{example}

\begin{remark}[Need for a measure solution theory]
\emph{The above example shows that the entropy solution concept
provided in Definition \ref{def:weak_sol} is only satisfactory for
short times, and it therefore needs to be improved to include also
measure solutions. On the other hand, allowing for a Dirac's delta
solution in \eqref{eq:main_oned} is complicated in view of the
nonlinearity w.r.t. $\rho$ in the drift term. In the next section
we shall solve this issue by producing a theory for measure
solutions. In particular, the concentration phenomenon (which is
somewhat hidden in the above computation) will be made evident.}
\end{remark}

\section{Measure solutions}\label{sec:measure}

The example provided in the example \ref{exe:explicit} shows that
a notion of measure solution is needed in order to establish a
global-in-time existence theory for \eqref{eq:main_oned}. In this
section we shall perform this task by considering the equation
satisfied by the pseudo-inverse of the distribution function of
the solution to \eqref{eq:main_oned}. Let us start by providing a
formal argument to justify the use of the pseudo inverse equation.
Assume $\rho$ to be the unique short--time smooth solution to
\eqref{eq:main_oned}, defined on $\R\times [0,T]$, with initial
datum $f_I$ having unit mass, and set
\begin{equation*}
    F(x,t):=\int_{-\infty}^x \rho(\xi,t)d\xi.
\end{equation*}
Let $[0,1]\times [0,T] \ni (z,t) \mapsto X(z,t)\in \R$ be the
pseudo-inverse of $F$
\begin{equation*}
  X(z,t)=\inf\{x:\ F(x,t)>z\}.
\end{equation*}
Formally, $X$ satisfies the equation
\begin{equation}\label{eq:pseudo_scaled}
    X_t = -\frac{X}{X_z^\gamma}.
\end{equation}
The interplay between \eqref{eq:pseudo_scaled} and
\eqref{eq:main_oned} will be more clear in the proof of Theorem
\ref{thm:main}. Clearly, the equation \eqref{eq:pseudo_scaled} in
the above form does not allow to consider a possible solution
$X(z,t)$ which is constant on some interval $I\subset [0,1]$. Let
us then instead consider the more general version
\begin{equation}\label{eq:pseudo_scaled_measures}
    X_t |X_z|^\gamma +X = 0,
\end{equation}
in which we have allowed for possible changes in the sign of
$X_z$. Please notice that having $X(\cdot,t)$ constant on some
interval $I\subset [0,1]$ is equivalent to a jump discontinuity
for $F(\cdot,t)$ at some point $x_0$, which means that $\rho
=\partial_x F$ is a probability measure with a non trivial Dirac
delta part at $x_0$. As already pointed out before in section
\ref{sec:multi_d}, we shall see that solutions to
\eqref{eq:main_oned} can only concentrate at $x=0$.

\begin{remark}[Connected components of $\rho$]\label{rem:connected}
\emph{Let us also mention that the equation
\eqref{eq:pseudo_scaled_measures} loses its meaning if the
solution $\rho(\cdot,t)$ to \eqref{eq:main_oned} is supported on a
disconnected set. Suppose for instance that the curve $x=x_1(t)$
is the left edge of a connected component of $\rho(\cdot,t)$ with
$x_1(t)<\sup\mathrm{supp}(\rho)$, and assume for simplicity that
no concentration occurs. Then, $X(\cdot,t)$ will be discontinuous
at the point $z_1=\int_{-\infty}^{x_1(t)}\rho(x,t)dx$ (which is
`locally' constant in time). In case of a finite number of
connected components, $X$ would be discontinuous on a set of
stationary points $z_1,\ldots,z_k$, until two connected components
possibly merge into one which means that $X$ becomes continuous on
one of the $z_j$'s. On the other hand, this creates unnecessary
complications to the notation and the calculations below.
Therefore, we shall restrict for simplicity to the case of an
initial condition $f_I$ with a \emph{single connected component in
its support}, for more details see Section \ref{sec:existence}.}
\end{remark}

Equation \eqref{eq:pseudo_scaled_measures} can be partially solved
via the method of characteristics for fully nonlinear first order
equations, see e. g. \cite{lions}. Let us recall that
\eqref{eq:pseudo_scaled_measures} is coupled with the initial
condition
\begin{equation*}
    X(z,0)=\xbar(z):=\inf\left\{x:\ \int_{-\infty}^x f_I(y)dy >z\right\},
\end{equation*}
where $z\in [0,1]$. We shall assume for simplicity that $\xbar \in \pseudospace \cap C^1((0,1))$.

Let us then define characteristics as follows: the vector
$Y(s):=(z(s), t(s))$ is the characteristic curve on the set
$[0,1]\times [0,T]$. We shall start by considering characteristics
generated at a point $(z_0,0)$ in the initial interval $t=0$. The
solution evaluated on the characteristics is $U(s):=
X(Y(s))=X(z(s),t(s))$. The gradient of $X$ along the
characteristics is given by
$P(s):=(p_1(s),p_2(s)):=(X_z(Y(s)),X_t(Y(s)))$. As for the initial
conditions at $s=0$, for a fixed initial point $z_0\in [0,1]$ we
have
\begin{align*}
    & Y(0)=(z_0,0),\qquad U(0)=\xbar (z_0)\\
    &
    P(0)=(p_1(0),p_2(0))=(\xbar_z(z_0),-\xbar(z_0)|\xbar_z(z_0)|^{-\gamma}).
\end{align*}
Let us differentiate \eqref{eq:pseudo_scaled_measures} with respect
to $z$ first and $t$ second:
\begin{align}
    & X_{tz}|X_z|^{\gamma} + \gamma \sign(X_z)X_t|X_z|^{\gamma-1} X_{zz} + X_z = 0 \label{eq:sec_der1}\\
    & X_{tt}|X_z|^{\gamma} + \gamma \sign(X_z)X_t|X_z|^{\gamma-1} X_{zt} + X_t = 0.\label{eq:sec_der2}
\end{align}
We have
\begin{align}
    & \dot{p}_1(s)= X_{zz}(Y(s))\dot{z}(s) + X_{zt}(Y(s))\dot{t}(s)\label{eq:der_p1}\\
    & \dot{p}_2(s)= X_{tz}(Y(s))\dot{z}(s) + X_{tt}(Y(s))\dot{t}(s).\label{eq:der_p2}
\end{align}
In order to remove the second derivatives in \eqref{eq:der_p1} and \eqref{eq:der_p2} we require
\begin{align}
    & \dot{z}(s) =\gamma \sign(X_z(Y(s)))X_t(Y(s))|X_z(Y(s))|^{\gamma-1}  = \gamma \sign(p_1(s))|p_1(s)|^{\gamma-1}p_2 (s),\label{eq:char_z}\\
    & \dot{t}(s) = |X_z (Y(s)))|^{\gamma} = |p_1 (s)|^{\gamma}\label{eq:char_t}
\end{align}
so that \eqref{eq:sec_der1} and \eqref{eq:sec_der2} imply
\begin{equation}\label{eq:char_p}
    \dot{p}_1(s) = - p_1(s),\qquad \dot{p}_2(s) = - p_2(s).
\end{equation}
Finally, differentiating $U(s)=X(z(s),t(s))$ with respect to $s$
yields
\begin{equation}\label{eq:char_u}
    \dot{U}(s) = p_1(s) \dot{z}(s) + p_2(s) \dot{t}(s) = (\gamma + 1) |p_1(s)|^\gamma p_2(s).
\end{equation}
Solving \eqref{eq:char_p} gives
\begin{equation*}
p_1(s)=\xbar_z(z_0)e^{-s}>0,\qquad
p_2(s)=-\xbar(z_0)\xbar_z(z_0)^{-\gamma}e^{-s},
\end{equation*}
and therefore \eqref{eq:char_u} gives $\dot U(s)=(\gamma +
1)\xbar(z_0) e^{-(1+\gamma) s}$, which yields
$U(s)=\xbar(z_0)e^{-(1+\gamma) s}$. Finally, solving
\eqref{eq:char_z} and \eqref{eq:char_t} yields
\begin{align*}
    & \dot{z}(s) =-\gamma \xbar(z_0)\xbar_z(z_0)^{-1}e^{-\gamma
    s},\\
    & \dot{t}(s) = \xbar_z(z_0)^{\gamma} e^{-\gamma s},
\end{align*}
which result in
\begin{align*}
    & z(s)=z_0 - \xbar(z_0)\xbar_z(z_0)^{-1}(1-e^{-\gamma s})\\
    & t(s)=\frac{1}{\gamma}\xbar_z(z_0)^{\gamma}(1-e^{-\gamma s}).
\end{align*}
We notice that $z$ can be written as a function of $t$ as follows:
\begin{equation*}
    z=z(t)=z_0-\gamma \xbar(z_0)\xbar_z(z_0)^{-(\gamma+1)} t
\end{equation*}

\begin{example}[Explicit solution revisited]\label{exe:example_rev}
\emph{Let us take as initial condition
\begin{equation*}
    \xbar (z)=z,
\end{equation*}
which corresponds (in terms of probability densities $\rho$) to
the initial condition in Section \ref{exe:explicit} up to a
dilation and a multiplication by a constant. The characteristic
curves obtained above are given by
\begin{align*}
    & U(s)=z_0 e^{-(1+\gamma)s}\\
    & z(s)=z_0 e^{-\gamma s}\\
    & t(s)=\frac{1-e^{-\gamma s}}{\gamma}\\
    & p_1(s)=e^{-s}\\
    & p_2(s)=-z_0 e^{-s}.
\end{align*}
Characteristics can be also written as
\begin{align*}
    & z=\tilde{z}(t)=z_0(1-\gamma t),\\
    & U=\tilde{U}(t)=z_0 (1-\gamma t)^{\frac{1+\gamma}{\gamma}}.
\end{align*}
The above characteristics occupy the region $A=\{(z,t)\in
[0,1]\times \R_+\ :\ 0\leq z \leq 1-\gamma t\}$. The formula for
the solution $X$ in $A$ is
\begin{equation}\label{eq:example_XonA}
  X(z,t)=z(1-\gamma t)^{1/\gamma}.
\end{equation}
In order to fill up the region $B:=[0,1]\times \R_+ \setminus A$,
we introduce the rarefaction wave originated at the point
$(z_0,t_0)=(1,0)$ in the next computations. \noindent Since
$U(0)=1$, assuming that the equation
\eqref{eq:pseudo_scaled_measures} is satisfied on the region $B$,
we have the following set of characteristics:
\begin{align}
   & U(s) = e^{-(\gamma + 1) s} \label{eq:example_rev1} \\
   & z(s) = 1+\lambda^{-1}(e^{-\gamma s}-1) \label{eq:example_rev2}\\
    & t(s) = \frac{\lambda^\gamma}{\gamma}(1-e^{-\gamma s})\label{eq:example_rev3}\\
    & p_1(s)=\lambda e^{-s} \nonumber\\
    & p_2(s)= -\lambda^{-\gamma}e^{-s},\nonumber
\end{align}
where $\lambda = p_1(0)$. Unlike in the region $A$, the value
$\lambda$ should not be prescribed according to the initial
condition, since the characteristic lines we are computing are
originated at $(1,0)$ and enter the region $B$. Therefore,
heuristically speaking we are not interested in the slope of the
initial condition at $z=1$, since the above mentioned
characteristic lines `will not touch' the profile originated from
the initial condition. In order to fill the whole region $B$, the
parameter $\lambda$ should be taken in the set $(1,+\infty)$. Such
parameter can be seen as the analogous of a self-similar variable.
\noindent In order to find the explicit formula for $X=X(z,t)$ on
$B$, we use \eqref{eq:example_rev2} and \eqref{eq:example_rev3} to
write $\lambda$ and $s$ in terms of $z$ and $t$. Inserting
\eqref{eq:example_rev3} in \eqref{eq:example_rev1} gives
\begin{equation}\label{eq:example_U1}
  U(s)=\widetilde{U}(t) =(1-\gamma \lambda^{-\gamma}t)_{+}^{\gamma + 1}{\gamma}.
\end{equation}
The positive part is needed in order to make formula
\eqref{eq:example_U1} well defined for all $t$.
\eqref{eq:example_rev2} provides the formula for $\lambda$, namely
$\lambda = \left(\frac{\gamma
t}{1-z}\right)^{\frac{1}{1+\gamma}}$, which can be inserted in
\eqref{eq:example_U1} to get to the formula for $X=X(z,t)$ on the
region $B$:
\begin{equation}\label{eq:exampleX}
  X(z,t)=\left[1-(\gamma t)^{\frac{1}{1+\gamma}}(1-z)^{\frac{\gamma}{1+\gamma}}\right]_{+}^{\frac{1+\gamma}{\gamma}}.
\end{equation}
From \eqref{eq:exampleX} and \eqref{eq:example_XonA} we deduce the following properties:
\begin{itemize}
  \item When $t<1/\gamma$, the two profiles for $X$ on $A$ and $B$ match in a $C^1$ way at the point $z=1-\gamma t$, with derivative equals to $\frac{\partial X}{\partial z}\big|_{z=1-\gamma t} = (1-\gamma t)^{1/\gamma}$.
  \item When $t\geq 1/\gamma$, only the profile for $X$ on $B$ survives. $X$ is constant zero on the interval $z\in [0,1-(\gamma t)^{-1/\gamma}]$, and $X(\cdot,t) \in C^1 [0,1)$ for all $t\geq 1/\gamma$.
  \item $X$ has an infinite $z$-derivative on $z=1$ for all times $t>0$.
\end{itemize}
Let $\mu(t)\in \mathcal{P}(\R)$ the distributional $x$-derivative
of $F(\cdot,t)$ the pseudo inverse of $X(\cdot,t)$. It is easily
seen that $\mu(t)=[1-(\gamma t)^{-1/\gamma}]\delta_0 +
\rho(\cdot,t)$ with $\rho(\cdot,t)\in BV(\R)$ being the weak
solution found in Example \ref{exe:explicit} except suitable time
change of variables to have unit mass. In particular, as
$t\rightarrow +\infty$, the concentrated mass of $\mu(t)$ tends to
the total mass $1$ as in the Example \ref{exe:explicit}. Moreover,
the concentrated mass is continuous in time.}
\end{example}

We now want to recover an equivalent formulation of the Oleinik
condition \eqref{eq:oleinik_simple} in terms of the pseudo inverse
variable $X$. Let $f_I\in \mespace$ and let $\rho$ be the
local-in-time solution to \eqref{eq:main_oned} provided by Theorem
\ref{thm:local} on $0\leq t< t^*$. Let $u_R$ solve \eqref{eq:IBV1}
and $u_L$ solve \eqref{eq:IBV2} respectively, and let $u$ be
defined as in \eqref{eq:u_def_split}. Notice that $\rho$ and $u$
are linked through the scaling \eqref{eq:change_to_burgers}. As
$\rho$ is bounded on $0\leq t< t^*$, then $u(0,t)=0$ for all
$0\leq t< t^*$. Denote
\begin{equation*}
  m_L=\int_{-\infty}^0 u(\xi,t) d\xi,\qquad m_R=1-m_L = \int_0^{+\infty} u(\xi,t) d\xi.
\end{equation*}
In view of $u(0,t)=0$, it is obvious that $m_L$ and $m_R$ are constant in time on $0\leq t < t^*$. Let us define the primitive variables
\begin{align*}
  & G_R(\xi,t)=\int_0^\xi u_L(\eta,t)d\eta,\quad \xi>0, \\
& G_L(\xi,t)=\int_{-\infty}^\xi u_R(\eta,t)d\eta \quad \xi<0,
\end{align*}
and their pseudo-inverses
\begin{align*}
   & Y_R(z,t) = \inf\{\xi>0:\ m_L + G_R(\xi,t)>z\},\quad z \in\ ]m_L,1], \\
   & Y_L(z,t) = \inf\{\xi<0:\ G_L(\xi,t)>z\},\quad z\in\ [0,m_L[, \\
& Y(z,t)=
\begin{cases}
Y_L(z,t) & z \in\ [0,m_L[\\
Y_R(z,t) & z\in\ [m_L,1]
\end{cases}.
\end{align*}
Since $\supp[u]$ is a connected interval and $u \in L^\infty(\R)$,
then both $Y_R$ and $Y_L$ are continuously differentiable in their
respective domains with the possible exceptions of the boundary
points $z=0$ and $z=1$ (for instance, $Y_R$ has an infinite
derivative at $z=1$ in case $\sup(\supp[u_R])>0$, since $u_R$
cannot have decreasing shocks on $\xi>0$, and it therefore tends
to zero at the right edge of its support).

Analogously to the $u$ variable, let us introduce left and right
distribution functions and pseudo-inverses for $\rho$. More
precisely, let
\begin{align}
  & F_R(x,t)=\int_0^x \rho(\eta,t)d\eta,\quad x>0, \label{eq:cumulativerhor}\\
& F_L(x,t)=\int_{-\infty}^x \rho(\eta,t)d\eta \quad x<0, \label{eq:cumulativerhol}\\
   & X_R(z,t) = \inf\{x:\ m_L + F_R(x,t)>z\},\quad z\in\ ]m_L,1],\nonumber \\
   & X_L(z,t) = \inf\{x:\ F_L(x,t)>z\},\quad z\in\ [0,m_L[,\nonumber\\
& X(z,t)=
\begin{cases}
X_L(z,t)& z\in\ [0,m_L[\\
X_R(z,t)& z\in\ ]m_L,1]
\end{cases}.\label{eq:pseudoXX}
\end{align}
By means of the scaling \eqref{eq:change_to_burgers}, a straightforward computation yields
\begin{equation}\label{eq:XtoY}
  Y(z,t)=\frac{\sign(X(z,t))}{\gamma}\left((1+\gamma)|X(z,t)|\right)^{\frac{\gamma}{1+\gamma}}.
\end{equation}
Since the boundary condition $u_L(0,t)=u_R(0,t)=0$ is satisfied, $Y$ solves the equation
\begin{equation*}
  Y_t Y_z^{\gamma}  +\sign(Y)\frac{1}{1+\gamma} = 0,\quad z\in\ [0,1], t>0.
\end{equation*}

We are now ready to rephrase the jump admissibility conditions
\eqref{eq:oleinik_simple} in terms of the pseudo inverse variable
$X$ defined in \eqref{eq:pseudoXX}. Assume for instance that
$X_z(\cdot,t)$ has an increasing jump on some point $z_0\in (0,1)$
with $X(z_0,t)>0$. Then, $Y(\cdot,t)$ defined in \eqref{eq:XtoY}
will have an increasing jump at the same point $z_0$. Formally,
since $G_R$ is the pseudo-inverse of $Y$ on the set $Y>0$, this
means that $G_R$ has a decreasing jump in its derivative, i. e.
$u_R$ has an decreasing jump, which is not admissible due to
\eqref{eq:oleinik_simple}. Thanks to this argument, and with a
specular one for the case $X(z_0,t)<0$, we can rephrase condition
\eqref{eq:oleinik_simple} as follows:
\begin{itemize}
  \item [(XJ)] Assume $X(\cdot,t)$ has a jump discontinuity in its $z$-derivative $X_z(\cdot,t)$ at some point $z_0$. Then
  \begin{align}
    & \lim_{z\nearrow z_0^-}X_z(z,t) < \lim_{z\searrow z_0^+}X_z(z,t)\qquad \hbox{if}\: X(z_0,t)<0,\nonumber\\
     &\lim_{z\nearrow z_0^-}X_z(z,t) > \lim_{z\searrow z_0^+}X_z(z,t)\qquad \hbox{if}\: X(z_0,t)>0,.
    \label{eq:oleinik_simple_rephrased}
\end{align}
\end{itemize}

\section{Existence and uniqueness of entropy measure solutions}\label{sec:existence}

We are now ready to state our notion of measure solution to \eqref{eq:pseudo_scaled_measures}. We shall use the following notation:
\begin{align*}
    \mespace = & \left\{\mu \in \mathcal{P}(\R):\ \mu = m\delta_0 + \rho\mathcal{L}_1,\; m\in [0,1],\, \rho \in L^1_+(\R),\, \supp(\rho) = [a,b],\; a,b\in \R\right\}.
\end{align*}
For a given $\mu \in \mespace$, with $\mu = m\delta_0 + \rho\mathcal{L}_1$ the pseudo-inverse distribution
\begin{equation}\label{eq:def_pseudo_measure}
  X_{\mu} (z)=\inf\{x\in \R:\ \mu((-\infty,x])>z\}
\end{equation}
satisfies $\mathrm{meas}(\{X_{\mu}(z)=0\}) = m$, and the following properties:
\begin{itemize}
  \item [(X1)] $X_{\mu}$ is continuous on $[0,1]$,
  \item [(X2)] $X_{\mu}$ is non-decreasing on $[0,1]$.
\end{itemize}
We shall use the notation
\begin{equation*}
  \pseudospace = \{X\in L^\infty([0,1]):\quad X\ \hbox{satisfies (X1)-(X2)}\}.
\end{equation*}
For future use, we also need the additional properties
\begin{itemize}
\item [(X3)] If $X_{\mu}(z)\neq 0$, then $X_{\mu}$ has finite
nonzero left and right derivatives $\partial_z^- X_{\mu}(z)$ and
$\partial_z^+ X_{\mu}(z)$.

\item [(X4)] If $X_\mu(0)<0$, then $\lim_{h\searrow 0} \frac{X_\mu
(h) - X_\mu(0)}{h} =+\infty$,

\item [(X5)] If $X_\mu(1)>0$, then $\lim_{h\searrow 0} \frac{X_\mu
(1) - X_\mu(1-h)}{h} =+\infty$,
\end{itemize}
and the notation
\begin{equation*}
  \pseudospace^+ = \{X\in \pseudospace:\quad X\ \hbox{satisfies (X3)--(X5)}\}.
\end{equation*}
Please notice that the conditions (X4) and (X5) are nothing but
avoiding non-entropic jump discontinuities at the edges of the
support of $\rho$. Condition (X3) will be used later
on to highlight the fact that $\rho$ is in
$BV_{loc}(\R\setminus\{0\})$.

\begin{defin}[Entropy measure solutions]\label{def:entropy_measure}
Let $f_I \in \mathcal{P}(\R)\cap BV(\R)$ with connected compact
support. A curve of probability measures $[0,+\infty)\ni t\mapsto
\mu(t)\in \mespace$ is an \emph{entropy measure solution} to
\eqref{eq:main_scaled} with initial datum $f_I$ if, given
$F(x,t):=\mu(t)((-\infty,x])$ its cumulative distribution, and the
pseudo-inverse function $X(\cdot,t):[0,1]\rightarrow \R$ of
$F(\cdot,t)$, the following properties are satisfied:
\begin{itemize}
\item [(In)] $\mu(0)=f_I\mathcal{L}_1$,

\item [(Re)] $X(\cdot,t) \in \pseudospace^+$ for all times $t\geq
0$,

\item [(Ol)] $X(\cdot,t)$ satisfies the rephrased Oleinik
condition \eqref{eq:oleinik_simple_rephrased} on $(0,1)$ and for
all $t>0$,

\item [(Eq)] $X(\cdot,t)$ satisfies
\eqref{eq:pseudo_scaled_measures} almost everywhere in $z\in
[0,1]$ and for all $t>0$.
\end{itemize}
\end{defin}

We now state the main result of our paper.

\begin{thm}[Existence and uniqueness of measures solutions]\label{thm:main}
Let $f_I \in \mathcal{P}(\R)\cap BV(\R)$ with connected compact
support. Then, there exists a unique global-in-time entropy
measure $[0,+\infty)\ni t\mapsto \mu(t)\in \mespace$ solution to
\eqref{eq:main_scaled} with initial datum $f_I$ in the sense of
Definition \ref{def:entropy_measure}. Moreover, the absolutely
continuous part $\rho \in L^\infty([0,\infty); L^1(\R))$ of
$\mu(t)$ is given by the change of variables
\eqref{eq:change_to_burgers} where $u$ is the unique entropy
solution to \eqref{eq:IBV1}--\eqref{eq:IBV2} in the sense of
Definition \ref{def:entropy_burgers}.
\end{thm}

\proof

\textsc{Step 1: construction of $u$.-} Let $\supp(f_I)=(A,B)$. For
simplicity we shall assume $a<0<b$, but the proof can be easily
repeated in the case of $a$ and $b$ having the same sign. Let
\begin{equation*}
  u_I(\xi):=(\gamma |\xi|)^{1/\gamma}f_I\left(\frac{\sign(\xi)}{1+\gamma}(\gamma
    |\xi|)^{(1+\gamma)/\gamma}\right).
\end{equation*}
Notice that $m_L:=\int_{-\infty}^0 u_I(\xi) d\xi  = \int_{-\infty}^0 f_I(x) dx$. Moreover, as $f_I \in L^\infty$, then $u_I$ is continuous at $\xi=0$ with $u_I(0)=0$. Let $u_L$ and $u_R$ be the unique (global in time) entropy solutions to \eqref{eq:IBV1} and \eqref{eq:IBV2} provided by Theorem \ref{thm:burgers} with initial conditions $u_{I,L}=u_I|_{\xi<0}$ and $u_{I,R}=u_I|_{\xi>0}$ respectively. Let $u$ be defined as in \eqref{eq:u_def_split}. By a simple comparison argument, in view of Remark \ref{rem:comparison}, we infer that $u(\xi)>0$ if and only if $\xi\in I(t)=(a(t),0)\cup (0,b(t))$ for some continuous functions $a(t)$ and $b(t)$. We set
\begin{equation*}
    t^*=\sup\{t\geq 0:\ u(\cdot,t)\: \hbox{is continuous at}\ \xi=0\ \hbox{and}\:\: u(0,t)=0\}.
\end{equation*}
Once again by a simple comparison argument, it is easily seen that $t^*>0$, since $u_I$ is continuous at $\xi=0$, and therefore the $u=0$ datum at $\xi=0$ is `transported' for at least some short time, until possible shocks may occur at $\xi=0$. Let
\begin{equation*}
  M_L(t):=\int_{-\infty}^0 u_L(\xi,t)d\xi,\qquad M_R(t):=\int_0^{+\infty} u_R(\xi,t)d\xi.
\end{equation*}
We remark that $M_L(t)$ and $M_R(t)$ are constant with $M_L(t) + M_R(t)=1$ on $0\leq t\leq t^*$.

\textsc{Step 2: construction of $\rho$.-} Now, let us set
\begin{align}
   & \rho_L(x,t)=\xi'(x) u_L(\xi(x),t),\quad \hbox{on}\ x<0, \label{eq:defrhol}\\
  & \rho_R(x,t)=\xi'(x) u_R(\xi(x),t),\quad \hbox{on}\ x>0, \label{eq:defrhor}\\
  & \rho(x,t)=
  \begin{cases}
  \rho_L(x,t) & \hbox{if}\ \ x<0\\
   \rho_R(x,t) & \hbox{if}\ \ x>0
  \end{cases},\label{eq:def_rho}
\end{align}
with $\xi(x)$ defined in \eqref{eq:change_x_to_xi}. Notice that
$\supp[\rho(\cdot,t)]\setminus \{0\}=\left(A(t),B(t)\right)
\setminus \{0\}$ for all times, with
\begin{align*}
    & A(t)= -\frac{1}{1+\gamma}(\gamma|a(t)|)^{\frac{1+\gamma}{\gamma}}\\
    & B(t)=\frac{1}{1+\gamma}(\gamma|b(t)|)^{\frac{1+\gamma}{\gamma}}.
\end{align*}
By comparison with the solution $\tilde{u}$ to
\eqref{eq:IBV1}--\eqref{eq:IBV2} with initial condition
$$
\tilde{u}_I=(\gamma
|\xi|)^{1/\gamma}\sup(f_I)\chi|_{(a(0),b(0))},
$$
from \eqref{eq:def_rho} we easily get
\begin{equation*}
    \rho(x,t)\leq C|x|^{-1/(1+\gamma)} \left(|x|^{\gamma/(1+\gamma)}\right)^{1/\gamma} =C
\end{equation*}
for $t$ small enough, which means that $\rho(\cdot,t)\in L^\infty$
for small times (see Example \ref{exe:explicit} for more details). On the other hand, since $\xi'(x)$ blows up at
$x=0$, $\rho$ may develop a blow-up at $x=0$ at time $t=t^*$ (see
Theorem \ref{thm:blow_up} below). However, we can estimate
\begin{equation}\label{eq:integrability_estimate}
    \rho(x,t) \leq |\xi'(x)|\|u(\cdot,t)\|_{L^\infty(\R)}\leq |\xi'(x)|\|u_I\|_{L^\infty(\R)} \leq C|x|^{-1/1+\gamma},
\end{equation}
where we have used $u_I \in L^\infty$, which is a trivial
consequence of the definition of $u_I$ and of $f_I$ having compact
support. As $\rho(\cdot,t)$ is compactly supported for all times,
then \eqref{eq:integrability_estimate} shows that
$\rho(\cdot,t)\in L^1(\R)$ for all $t\geq 0$. Moreover,
$x\rho^{1+\gamma}\in L^\infty(\R)$. We claim that $\rho$ is a weak
solution to \eqref{eq:main_oned} in the sense if Definition
\ref{def:weak_sol}. To see this, we have to verify that $\rho$
satisfies \eqref{eq:def_weak} for all $\phi\in
C_c^\infty(\R\times[0,+\infty))$ with $\phi(0,t)=0$. The first and
the second term in \eqref{eq:def_weak} scale in a straightforward
way. Let us compute, for $\widetilde{\phi}(\xi,t)=\phi(x(\xi),t)$,
\begin{align*}
    & \int_0^{+\infty}\int_0^{+\infty} \rho^{1+\gamma}(x,t) x \phi_x(x,t) dx dt = \int_0^{+\infty}\int_0^{+\infty} \xi'(x)^{1+\gamma} x u_R^{1+\gamma}(\xi,t)\phi_x(x,t) dx dt \\
    & \ = \int_0^{+\infty}\int_0^{+\infty} [(1+\gamma)x]^{-1} x u_R^{1+\gamma}(\xi(x),t)\xi'(x)\widetilde{\phi}_\xi(\xi(x),t)dx dt\\
     & \ = \int_0^{+\infty}\int_0^{+\infty} (1+\gamma)^{-1} u_R^{1+\gamma}(\xi,t)\widetilde{\phi}_\xi(\xi,t)d\xi dt,
\end{align*}
and this term, combined with the other two terms and with the
definition of weak solution for $u$, and with the terms obtained
on $x<0$, gives the desired formula \eqref{eq:def_weak}. Since $u$
is $BV(\R)$ (cf. \cite{dafermos}), the right and left limits of $u$
are always defined, and $u$ has at most a countable number of
jumps. The jumps are decreasing on $\xi<0$ and increasing on
$\xi>0$ since $u$ is an entropy solution. The scaling
$u\rightarrow \rho$ clearly preserves the ordering of the jumps,
and hence all the jumps of $\rho$ are admissible. We stress that
the $\rho$ introduced here is globally defined for all times. Such
$\rho$ represents the absolutely continuous part of the candidate
measure solution to our problem.

\textsc{Step 3: construction of the cumulative distributions.-} As
$\rho\in L^1(\R)$, we can define the cumulative distributions
$F_L$ and $F_R$ as in \eqref{eq:cumulativerhol} and
\eqref{eq:cumulativerhor}. Please notice that $F_L$ (resp. $F_R$)
is a strictly increasing bijection from $[A(t),0)$ (resp.
$(0,B(t)]$) onto $[0,M_L(t))$ (resp. $(0,M_R(t)]$). Moreover,
$F_R$ and $F_L$ are defined globally in time. Let $0<t_1<t_2$ and
let us apply Definition \ref{def:weak_sol} with the test functions
$\varphi_1^\epsilon(x,t)=\chi^\epsilon_{[t_1,t_2]}(t)\chi^\epsilon_{(-\infty,x_0]}(x)$
for $x_0<0$ and
$\varphi_2^\epsilon(x,t)=\chi^\epsilon_{[t_1,t_2]}(t)\chi^\epsilon_{[x_0,+\infty)}(x)$
for $x_0>0$, where $\chi^\epsilon$ is a standard mollification of
a characteristic function. After sending $\epsilon \rightarrow 0$,
multiplying by $(t_2-t_1)^{-1}$, and sending $t_2\searrow t_1$ we
get
\begin{align}
    & \frac{\partial F_L(x_0,t_1)}{\partial t} = x_0 \rho_L(x_0,t_1)^{1+\gamma}\quad x_0<0,\label{eq:weak_cumulative1}\\
    & \frac{\partial F_R(x_0,t_1)}{\partial t}-M_R'(t) = x_0 \rho_R(x_0,t_1)^{1+\gamma}\quad x_0>0.\label{eq:weak_cumulative2}
\end{align}

\textsc{Step 4: definition of the candidate solution $X(z,t)$.-}
Let us set
\begin{align*}
  & X_L(z,t)=\inf\{x:\ F_L(x,t)>z\}\qquad \hbox{on}\ 0\leq z < M_L(t) \\
  & X_R(z,t)=\inf\{x:\ 1-M_R(t) + F_R(x,t)>z\} \qquad \hbox{on}\ 1-M_R(t) < z < 1,
\end{align*}
and
\begin{equation}\label{eq:pseudo_solution}
  X(z,t)=
\begin{cases}
X_L(z,t) & \hbox{on}\ 0\leq z < M_L(t)\\
0 & \hbox{on}\ M_L(t)\leq z \leq 1-M_R(t)\\
X_R(z,t) & \hbox{on}\ 1-M_R(t) < z < 1
\end{cases}.
\end{equation}
The \emph{candidate solution} $\mu(t)\in \mathcal{M}_0$ is defined
as $\mu(t)=\partial_x F(\cdot,t)$ in the sense of distributions,
with $F(x,t)=\inf\{z\in [0,1]\ :\ X(z,t)>x\}$ being the pseudo
inverse of $X$. Notice in particular that $X=X_\mu$ according to
\eqref{eq:def_pseudo_measure}. It is clear from the above
definition that $\mu(t)\in \mathcal{M}_0$, with $\mu(t) =
m(t)\delta_0 + \rho(\cdot,t)\mathcal{L}_1$, which proves that the
absolutely continuous part of $\mu(t)$ is $\rho(t)$ given in Step
3.

As $f_I\in \mespace$, then $X(\cdot,0)\in \pseudospace$, and the property (In) is trivially proven. Now, as $X_L(\cdot,t)$ and $X_R(\cdot,t)$ are both strictly increasing, we have that $X_L(\cdot,t):[0,M_L(t))\rightarrow [a(t),0)$ is the inverse of $F_L(\cdot,t)$ and $X_R(\cdot,t):(1-M_R(t),1]\rightarrow (0,b(t)]$ is the inverse of $1-M_R(t)+F_R(\cdot,t)$. Therefore,
\begin{align}
    & \partial_z X_L(z,t)= (\rho_L(X_L(z,t),t))^{-1},\quad z\in (0,M_L(t)),\label{eq:prepseudo1}\\
    & \partial_z X_R(z,t)= (\rho_R(X_R(z,t),t))^{-1},\quad z\in (1-M_R(t),1).\label{eq:prepseudo2}
\end{align}
Moreover
\begin{align}
    & \partial_t X_L(z,t)= -\partial_z X_L(z,t) \partial_t F_L(X(z,t)t),\quad z\in (0,M_L(t)),\label{eq:prepseudo3}\\
    & \partial_t X_R(z,t)= -\partial_z X_L(z,t) \partial_t \left( F_R(X(z,t),t) - M_R(t)\right),\quad z\in (1-M_R(t), 1).\label{eq:prepseudo4}
\end{align}
We can then combine \eqref{eq:weak_cumulative1}--\eqref{eq:weak_cumulative2} with \eqref{eq:prepseudo1}--\eqref{eq:prepseudo4} to obtain that \eqref{eq:pseudo_scaled_measures} is satisfied on the set $z\in (0,M_L(t))\cup (1-M_R(t),1)$. Since $X$ defined in \eqref{eq:pseudo_solution} is constantly zero in the remaining set $z\in [M_L(t),1-M_R(t)]$, then the property (Eq) is satisfied.

\textsc{Step 5: regularity.-} As for the property (Re), it is
clear that $X(\cdot,t)$ is continuous on $z\notin
[M_L(t),1-M_R(t)]$. Moreover, since $F_L(0,t)=M_L(t)$ and
$F_R(0,t)=0$, then $X_L(M_L(t)^-,t)=0=X_R((1-M_R(t))^+,t)$, which
implies continuity on the whole set $[0,1]$. Moreover,
$X(\cdot,t)$ is bounded, and non-increasing. Since $u(\cdot,t)\in
BV(\R)$ for all times, then $\rho(\cdot,t)\in
BV_{loc}(\R\setminus\{0\})$, and therefore the property (X3) is
proven. Since $u$ is an entropy solution in the sense of
Definition \ref{def:entropy_burgers}, an increasing jump at
$\xi=a(t)$ is not entropic, and therefore the same holds for
$\rho(\cdot,t)$ at $x=A(t)$. Therefore, $\rho(\cdot,t)$ is
continuous at $x=A(t)$ with $\rho(A(t),t)=0$. Similarly, one can
prove that $\rho(B(t),t)=0$. This proves properties (X4) and (X5).
Thus, (Re) is proven.

Finally, we prove the property (Ol), namely that the Oleinik
condition \eqref{eq:oleinik_simple_rephrased} is satisfied. Let
$z_0$ a point of jump for $X_z(\cdot,t)$, and let $x_0=X(z_0,t)$.
Assume first that $x_0<0$. It is clear that $\rho_L$ has a jump at
$x_0$. Since $\rho$ is entropic, the jump is decreasing. Since
$X_z(z_0^-,t)=(\rho_L(x_0^-,t)^{-1}$ has an increasing jump, as
required in \eqref{eq:oleinik_simple_rephrased}. A similar
computation holds in the case $x_0>0$.

\textsc{Step 6: uniqueness.-} Finally, we prove uniqueness. Let
$\Xtil$ be another solution in the sense of Definition
\ref{def:entropy_measure} with initial condition $f_I$. Let
\begin{align*}
    & I_L(t)=\{z\in [0,1]\ :\; \Xtil(z,t)<0\},\\
    & I_R(t)=\{z\in [0,1]\ :\; \Xtil(z,t)>0\}.
\end{align*}
Let
\begin{align*}
    \Ftil_L(x,t)=\inf \{z\in I_L(t)\ :\; \Xtil(z,t)>x\},\qquad \hbox{on}\ \ x<0\\
    \Ftil_R(x,t)=\inf \{z\in I_R(t)\ :\; \Xtil(z,t)>x\},\qquad \hbox{on}\ \ x>0.
\end{align*}
Both $\Ftil_L$ and $\Ftil_R$ are strictly increasing on their
respective domains because of property (X1). As a consequence of
\eqref{eq:pseudo_scaled_measures}, a simple computation shows that
\begin{align*}
    & \partial_t\Ftil_{L}(x,t)=(\partial_x\Ftil_{L})^{1+\gamma} x,
\end{align*}
for all $x<0$. Therefore, $\rhotil_L:=\Ftil_{L,x}$ satisfies
\eqref{eq:main_oned} on $x<0$, and a symmetric argument shows that
$\rhotil_R:=\Ftil_{R,x}$ satisfies \eqref{eq:main_oned} on $x>0$.
Let
\begin{equation*}
    \rhotil(x,t)=
    \begin{cases}
    \rhotil_L(x,t) &\hbox{if}\ \ x<0\\
    \rhotil_R(x,t) &\hbox{if}\ \ x>0
    \end{cases},
\end{equation*}
and let $\util$ be defined by
$\util(\xi,t)=x(\xi)\rhotil(x(\xi),t)$ with $x(\xi)$ given as in
\eqref{eq:change_inverse}. The usual scaling computation gives
that $\util$ is a weak entropy solution to
\eqref{eq:IBV1}--\eqref{eq:IBV2} in the sense of Definition
\ref{def:entropy_burgers}. The fact that $\util$ satisfies the
entropy condition \eqref{eq:oleinik_simple} comes from property
\eqref{eq:oleinik_simple_rephrased} and from the reversed argument
at the end of Step 5 above. Therefore, as $\Xtil(z,0)=X(z,0)$
implies that $u$ and $\util$ have the same initial condition, the
result in Theorem \ref{thm:burgers} implies that $u\equiv \util$.
This also implies that $I_L(t)=[0,M_L(t))$ and
$I_R(t)=(1-M_R(t),1]$, which imply also that $X$ and $\Xtil$
coincide on the set $z\neq 0$. Now, at each time $t$ there are two
possible situations: either $I_L(t)\cup I_R(t)=[0,1]$, or
$I_L(t)\cup I_R(t)\subsetneqq[0,1]$. In the former case,
$\Xtil(\cdot,t)\equiv X(\cdot,t)$ on $[0,1]$. In the latter case,
as we clearly have $\Xtil(z,t)\nearrow 0$ as $z\nearrow \sup
I_L(t)$ and $\Xtil(z,t)\searrow 0$ as $z\searrow \inf I_R(t)$, the
only possible way to extend $\Xtil$ on the set $[0,1]\setminus
(I_L(t)\cup I_R(t))$ in a way such that $\Xtil(\cdot,t)$ is
non-decreasing is to set $\Xtil(z,t)= 0$ on $z\in [0,1]\setminus
(I_L(t)\cup I_R(t))$. Therefore, $\Xtil$ coincides with $X$
defined in \eqref{eq:pseudo_solution}, and the proof is complete.
\endproof

\begin{remark}
\emph{In theorem \ref{thm:main} we have assumed for simplicity
that the support of the initial condition has just one connected
component. Such assumption could be easily removed by considering
a finite number of components, and then generalizing the above
result to all $BV$ initial probability densities with compact
support by approximation. We omit the details.}
\end{remark}

\begin{remark}[Existence and uniqueness for \eqref{eq:main}]\label{rem:original1}
\emph{The definition \ref{def:entropy_measure} and the statements
in Theorem \ref{thm:main} can be easily transferred to the level
of $f$ solution to \eqref{eq:main}. Roughly speaking, a notion of
measure solution to \eqref{eq:main} can be recovered, which
satisfies the same properties of $\mu(t)$ in Theorem
\ref{thm:main}.}
\end{remark}

\section{Qualitative and asymptotic behaviour}\label{sec:qualitative}

We now prove that all measure solutions to \eqref{eq:main_oned}
develop a singularity in finite time, and that the concentrated
mass is strictly increasing in time, and it converges to the total
mass for large times.

\begin{thm}[Finite time blow up]\label{thm:blow_up}
Let $f_I \in \mathcal{P}(\R)\cap BV(\R)$ with connected compact
support. Then, the unique global-in-time entropy measure solution
$\mu(t) = m(t)\delta_0 + \rho(\cdot,t)\mathcal{L}_1$ to
\eqref{eq:main_scaled} with initial datum $f_I$ in the sense of
Definition \ref{def:entropy_measure} satisfies the following
properties:
\begin{enumerate}
  \item [(FTBU)] There exists a time $T<+\infty$ such that $\rho(\cdot,t)\not \in L^\infty(\R)$ for all $t\geq T$.
  \item [(SCon)] The function $m(\cdot)$ is continuous and strictly increasing in time, and it satisfies $$\lim_{t\rightarrow +\infty}m(t)=1.$$
\end{enumerate}
\end{thm}

\proof The proof relies on an extensive use of the comparison
principle \ref{rem:comparison}. Let $\supp(f_I)=(a,b)$, and assume
for simplicity that $b>0$. Let $[c,d]\subset (0,b)$ such that
\begin{equation*}
    \lambda:=\min_{x\in[c,d]} f_I(x)>0.
\end{equation*}
We define $u_I(\xi)$ as in the proof of Theorem \ref{thm:main},
\begin{equation*}
  u_I(\xi):=(\gamma |\xi|)^{1/\gamma}f_I\left(\frac{\sign(\xi)}{1+\gamma}(\gamma
    |\xi|)^{(1+\gamma)/\gamma}\right),
\end{equation*}
and consider $u$ the unique entropy solution to \eqref{eq:IBV1}--\eqref{eq:IBV2} with initial condition $u_I$. Let
\begin{equation*}
    \bar u_0(\xi):=\lambda(\gamma |\xi|)^{1/\gamma}\chi_{[c,d]},
\end{equation*}
and consider the corresponding entropy solution $\bar u_R$ to
\eqref{eq:IBV1} on $\xi>0$. We next show that $\bar u(0^+,t)$
becomes nonzero in a finite time $t^*$. The solution $\bar u$ has
an increasing shock wave $\xi=s(t)$ originating at $\xi=c$, with
velocity $-\bar u(s(t))^\gamma/(1+\gamma)$. At $\xi=d$, a
rarefaction wave is originated with left front given by the
straight line $\xi(t)= d - t \lambda^\gamma \gamma d$, i. e. the
rarefaction wave is faster than the shock wave. If the shock wave
reaches $\xi=0$ before the rarefaction wave does, then the
assertion is proven. Otherwise, another shock wave is originated,
which can be easily computed to travel with a law of order
$-Ct^{1/\gamma}$ for large times. Hence, the rarefaction wave
reaches zero in a finite time, and the assertion that $\bar
u(0^+,t)$ becomes nonzero in a finite time $t^*$ is proven. Now,
we employ the comparison principle in Remark \ref{rem:comparison}
to deduce that $u(0^+,t)$ also becomes nonzero at least for $t\geq
t^*$. Now, if $\mu(t)$ has already developed a blow-up singularity
in $L^\infty$ before $t=t^*$, then (FTBU) is trivially true.
Assume this is not the case, then we know from the proof of
Theorem \ref{thm:main} that $\rho$ defined in \eqref{eq:def_rho}
is the unique entropy solution to \eqref{eq:main_oned}, and the
above computations show that $\lim_{x\searrow
0^+}\rho(x,t^*)=+\infty$, which proves (FTBU).

Now, a trivial application of the definition of weak solution for
$u$ (with suitable choice of the test functions and with suitable
mollifiers) gives, for all $s>t$
\begin{equation}\label{eq:mass_flux}
    \int_0^{+\infty} u(\xi,s) d\xi - \int_0^{+\infty} u(\xi,t) d\xi = -\int_t^s u^{1+\gamma}(0^+,\sigma)d\sigma.
\end{equation}
Now, let $t_1=\sup\{t>0\ :\ \ m(t)=0\}$. The support of
$u(\cdot,t)$ is a connected interval, and we know that
$u(0^+,t_1)>0$. To see this, assume by contradiction that this is
not the case, then there exists an interval $(0,\delta)$ on which
$u(\cdot,t)$ is zero, and by a simple characteristics argument it
is easy to see that $u(0^+,t_1)$ will remain zero for some short
time $[t_1,t_1+\epsilon)$ with $\epsilon>0$, which contradicts the
maximality of $t_1$. Therefore, $\inf\supp(u_R(\cdot,t))=0$ for
all $t\geq t_1$. Once again, a very simple characteristics
argument shows that $u(0^+,t)>0$ for all $t \geq t_1$, as in that
point $u$ will achieve values originated from the interior of
$\supp(u(\cdot,t_1))$. Hence, \eqref{eq:mass_flux} implies that
\begin{equation*}
     \int_0^{a(s)} u(\xi,s) d\xi < \int_0^{a(t)} u(\xi,t) d\xi,
\end{equation*}
which gives
\begin{equation*}
     \int_0^{+\infty} \rho(x,s) dx < \int_0^{+\infty} \rho(x,t) dx,
\end{equation*}
and hence $m(t)$ is increasing by conservation of the total mass.
To finish the proof of (SCon), we need to prove that all the mass
concentrates as $t\rightarrow+\infty$. To prove this assertion,
let $C=\|f_I\|_{L^\infty}$, and let $[A,B]\supset \supp(f_I)$.
Define
\begin{equation*}
    \bar{f}(x)= C\chi_{[A,B]},
\end{equation*}
and consider
\begin{equation*}
  \widetilde{u}_I(\xi):=(\gamma |\xi|)^{1/\gamma}\bar{f}\left(\frac{\sign(\xi)}{1+\gamma}(\gamma
    |\xi|)^{(1+\gamma)/\gamma}\right).
\end{equation*}
Let $\widetilde{u}(\xi,t)$ be the entropy solution in the sense of
Definition \ref{def:entropy_burgers} with initial condition
$\widetilde{u}$. A simple characteristic argument shows that
\begin{equation*}
    \lim_{t\rightarrow +\infty} \|\widetilde{u}(\cdot,t) \|_{L^\infty} = 0.
\end{equation*}
Hence, by the usual comparison argument,
\begin{equation*}
    \lim_{t\rightarrow +\infty} \|u(\cdot,t) \|_{L^\infty} = 0.
\end{equation*}
Still by a comparison argument, one easily sees that
$\supp(u(\cdot,t))$ is uniformly bounded in time, and this implies
\begin{equation*}
  \int_\R u(\xi,t) d\xi \rightarrow 0,\qquad \hbox{as}\quad t\rightarrow +\infty,
\end{equation*}
From the proof of Theorem \ref{thm:main}, $\rho$ defined in
\eqref{eq:def_rho} from $u$ is the absolutely continuous part of
$\mu(t)$, and we have
\begin{equation*}
  \int_\R\rho(x,t) dx = \int_\R u(\xi,t) d\xi \rightarrow 0,\qquad \hbox{as}\quad t\rightarrow +\infty,
\end{equation*}
which implies the assertion.
\endproof

As a straightforward consequence of the above proof, we have the
following

\begin{corollary}\label{cor:support}
Under the same assumptions of Theorem \ref{thm:blow_up}, the
support of $\mu(t)$ is uniformly bounded in time.
\end{corollary}

\begin{remark}
\emph{The strict monotonicity of $m(t)$ would be clearly violated
in case of solutions with many components in their support. In
that case, the concentrated mass would stay constant for some
waiting time until a new wave will hit the boundary $x=0$.}
\end{remark}

\begin{remark}[Asymptotic behaviour for the original equation \eqref{eq:main}]\label{rem:original2}
\emph{A similar statement to that in Theorem \ref{thm:blow_up} can
be formulated for the global measure solution to \eqref{eq:main}
mentioned in Remark \ref{rem:original1}. Once again, the same
qualitative properties can be easily recovered via the scaling
\eqref{eq:scaling}. Additionally, the measure of the support of
all measure solutions of \eqref{eq:main} converges algebraically
to zero for large times. As a consequence of the scaling
\eqref{eq:scaling} and of Corollary \ref{cor:support}, we can
easily see that
\begin{equation*}
    W_p(\delta_0,f(t))\leq C (1+t)^{-1/\gamma},
\end{equation*}
for some $C>0$, and for all $1\leq p\leq \infty$. Here $W_p$
denotes the $p$-Wasserstein distance, see e. g. \cite{villani}.
$f(t)$ is, by abuse of notation, the global-in-time measure
solution to \eqref{eq:main} mentioned in Remark
\ref{rem:original1}.}
\end{remark}

\end{document}